\documentclass{amsart} 

\usepackage{amsfonts}
\usepackage{amsmath}
\usepackage{amsthm}

\newcommand{\R}{\operatorname{Re}}

\newcommand{\la}{\lambda}

\newcommand{\ar}{\right|}
\newcommand{\al}{\left|}

\newtheorem{twot}{Lemma}
\newtheorem{cmvt}[twot]{Theorem}
\newtheorem{reim}[twot]{Lemma}

\newtheorem{enequals2}[twot]{Lemma}

\newtheorem{vd1}[twot]{Corollary}
\newtheorem{vd}[twot]{Lemma}
\newtheorem{cor}[twot]{Corollary}
\theoremstyle{plain}

\newtheorem{mvt}[twot]{Lemma}

\newtheorem{divdif}[twot]{Theorem}
\usepackage{amscd}
\usepackage{amscd}

\begin{document}
\address{School of Mathematics,
University of New South Wales,
Sydney, NSW 2052,
Australia}

\email{keith@maths.unsw.edu.au}

\author{Keith Rogers}
\date{}
\title{Sharp van der Corput estimates and minimal divided differences}

\begin{abstract}We find the nodes that minimise divided differences and use them to find the sharp constant in a sublevel set estimate. We also find the sharp constant in the first
instance of the van der Corput Lemma using a complex mean value theorem for integrals. With these sharp bounds we improve the constant in the general van
der Corput Lemma, so that it is asymptotically sharp.
\end{abstract}
\maketitle

\section{Sublevel set estimates}

The importance of sublevel set estimates in van der Corput lemmas was highlighted by A. Carbery, M. Christ, and J. Wright \cite{ca}, \cite{cw}. We find the sharp constant in the
 following sublevel set estimate. We will use this in Section~\ref{threet} to prove an asymptotically sharp van der Corput lemma. The constants $C_n$ take different values in
each lemma. 

\begin{twot}\label{twot} Suppose that $f:(a,b)\to \mathbb{R}$ is $n$ times differentiable with $n\ge1$ and $|f^{(n)}(x)|\ge \la >0$
on $(a,b).$ Then
$$
|\{ x \in (a,b) : |f(x)|\le \alpha\}| \le C_n\left(\alpha/\lambda\right)^{1/n},
$$
where $C_n=(n!2^{2n-1})^{1/n}.$
\end{twot}
\noindent We note that $C_n\le 2n$ for all $n\ge 1,$ and by Stirling's formula,
$$
\lim_{n\to \infty} C_n-4n/e=0.
$$


The Chebyshev polynomials will be key to the proof of Lemma~\ref{twot}, so we recall some facts that we will need. For a more complete introduction see \cite{th}.
Consider $T_n$ defined on $[-1,1]$ by
$$
T_n(\cos\theta)=\cos n \theta,
$$
where $ 0\le\theta \le \pi$. If we take the binomial expansion of
$$
\cos n\theta +i\sin n \theta=\left(\cos \theta +i\sin \theta\right)^n
$$
and equate real parts, we obtain
$$
\cos n\theta=\sum_{k=0}^{\lfloor n/2\rfloor}\left((-1)^k\sum_{j=k}^{\lfloor n/2\rfloor} {n \choose 2j }{j \choose k}\right)(\cos\theta)^{n-2k},
$$
where $\lfloor n/2\rfloor$ denotes the integer part of $n/2.$ Thus we can consider $T_n$ to be the polynomial of degree $n$ defined on the real line by
$$
T_n(x)=\sum_{k=0}^{\lfloor n/2\rfloor}\left((-1)^k\sum_{j=k}^{\lfloor n/2\rfloor} {n \choose 2j }{j \choose k}\right)x^{n-2k}.
$$
It is clear that $|T_n|\le 1$ on $[-1,1],$ and has extrema at $\eta_j=\cos j\pi/n$ for
$j=0,\ldots,n.$ Finally we calculate the leading coefficient,
$$
\sum_{j=0}^{\lfloor n/2\rfloor}{n \choose 2j }=\frac{1}{2}\left((1+1)^n+(1-1)^n\right)=2^{n-1}.
$$

We will also require the following generalisation of
the classical mean value
theorem. It can be found in texts on numerical analysis, for example \cite{ik}. We include a proof for convenience.

\begin{mvt}\label{mvt} Suppose that $f:(a,b)\to \mathbb{R}$ is $n$ times differentiable, where $n\ge1,$ and suppose 
that $x_0<x_1<\ldots<x_{n}$ are distinct points in $(a,b).$
Then there exists $\zeta\in (a,b)$ such that
$$
f^{(n)}(\zeta)=\sum_{j=0}^{n}c_jf(x_j),
$$
where $
c_j = 
(-1)^{j+n} n! \prod_{k:k\neq j} |x_k-x_j|^{-1}.$
\end{mvt} 
\begin{proof} Let $\phi (x) = f(x) -\mu x^n-p(x)$ where $\mu\in \mathbb{R}$ and $p$ is a polynomial of degree $n-1.$   If we can choose $\mu$ and $p$ such that
$\phi(x_j)=0$ for $j=0,\ldots ,n,$ then $\phi^{(n)}$ must vanish somewhere in $(a,b),$ by Rolle's Theorem.  Thus, there will exist $\zeta\in (a,b)$ such that
$f^{(n)}(\zeta)=n!\mu.$ 

It remains to choose $\mu$ and $p$ so that $\phi(x_j)=0$ for $j=0,\ldots,n.$ We require a solution to
$$
V_{n+1}(x_0,\ldots,x_n)\left(\begin{array}{c}
a_0\\
\vdots\\
a_{n-1}\\
\mu\end{array}\right)=\left(\begin{array}{c}
f(x_0)\\
\vdots\\
\vdots\\
f(x_n)\end{array}\right),
$$
where $a_0,\ldots,a_{n-1}$ are the coefficients of $p$ and $V_{n+1}$ is the Vandermonde matrix, defined by
$$
V_{n+1}(x_0,\ldots,x_n)=\left(\begin{array}{ccccc}
1&x_0&x_0^2&\ldots&x_0^n\\
1&x_1&x_1^2&\ldots&x_1^n\\
\vdots&\vdots&\vdots& &\vdots\\
1&x_n&x_n^2&\ldots&x_n^n
\end{array}\right).
$$

There is clearly a unique solution when the $x_j$ are distinct, and by Cram\'{e}r's rule, we have
$$
\mu=
\sum_{j=0}^n(-1)^{j+n} \frac{\det V_n(x_0,\ldots,x_{j-1},x_{j+1},\ldots,x_n)}{\det V_{n+1}(x_0,\ldots,x_n)}f(x_j).
$$
 Finally, the determinant of the Vandermonde matrix is given by
$$
\det V_{m+1}(y_0,\ldots,y_m)=\prod_{0\le k<j\le m}(y_k-y_j),
$$
so that $\mu$ is as desired.
\end{proof}

We apply this result to the Chebyshev polynomials. It is easy to see that $T_n^{(n)}=n!2^{n-1},$ and that  $
T_n(\eta_j)=
(-1)^{j+n}$ for the Chebyshev extrema $\eta_j=\cos j\pi/n$.
Thus we obtain
\begin{equation*}
\sum_{j=0}^n\prod_{k:k\neq j} n!|\eta_k-\eta_j|^{-1}= n!2^{n-1},
\end{equation*}
 so that
\begin{equation}\label{po}
\sum_{j=0}^n\prod_{k:k\neq j}|\eta_k-\eta_j|^{-1}= 2^{n-1}.
\end{equation}
\ \newline

\noindent\textit{Proof of Lemma \ref{twot}.}
Suppose that $|E|=|\{ x \in (a,b) : f(x)\le \alpha\}|>0.$ First we map $E$ to an interval of the same measure without increasing distance. Centering at the origin and scaling by $2/|E|$ if
necessary, we map the interval to $[-1,1].$
Now by (\ref{po}) we have
$$
\sum_{j=0}^n\prod_{k:k\neq j} |\eta_k-\eta_j|^{-1}= 2^{n-1},
$$
where $\eta_j=\cos j\pi/n.$ Thus, mapping back to $E,$ 
there exist $n+1$ points $x_0,\ldots,x_n\in E$ such that
$$
\sum_{j=0}^n\prod_{k:k\neq j} |x_k-x_j|^{-1}\le 2^{n-1}\frac{2^{n}}{|E|^n}=\frac{2^{2n-1}}{|E|^n}.
$$
On the other hand, by Lemma~\ref{mvt} there exists $\zeta\in (a,b)$ such that
$$
f^{(n)}(\zeta)=\sum_{j=0}^{n}c_jf(x_j),
$$
where $|c_j| =  n! \prod_{k:k\neq j} |x_k-x_j|^{-1}.$ Putting these together, 
\begin{equation*}
\lambda\le \left|\sum_{j=0}^{n}c_jf(x_j)\right|\le \sum_{j=0}^{n}n!\prod_{k:k\neq j} |x_k-x_j|^{-1}\alpha\le n!\frac{2^{2n-1}}{|E|^n}\alpha,
\end{equation*}
so that 
$$
|E|\le \left(n!2^{2n-1}\right)^{1/n}\left(\alpha/\lambda\right)^{1/n},
$$
as desired.\hfill $\square$

The sharpness may be observed by considering $f(x)=T_n(x)$ on $[-1,1]$ with $\la =n!2^{n-1}$ and $\alpha=1.$

\section{Divided differences}

Divided differences were first considered by Newton and are important in interpolation theory. For a given
set of nodes $x_0<\ldots<x_m\in[-1,1],$ the $n$th divided difference is defined
recursively by
$$
f[x_0,\ldots,x_n]=\frac{f[x_0,\ldots,x_{n-1}]-f[x_1,\ldots,x_{n}]}{x_0-x_n},
$$ 
where $f[x_0]=f(x_0)$ and $f[x_1]=f(x_1).$ It is not difficult to calculate that
$$
f[x_0,\ldots,x_n]= 
\sum_{j=0}^n(-1)^{j+n}\prod_{k:k\neq j} |x_k-x_j|^{-1}f(x_j).
$$

The following shows that the Chebyshev extrema are optimal for the minimisation of divided differences.
\begin{divdif} The Chebyshev extrema, $\eta_j=\cos j\pi/n$ for $j=0,\ldots,n,$ are the unique nodes where 
$$
f[\eta_0,\ldots,\eta_n]\le 2^{n-1}||f||_{\infty},
$$
for all $f\in C[-1,1].$ 
\end{divdif}

\begin{proof} That the inequality holds is clear from (\ref{po}).  To show the uniqueness, we suppose there exist $x_0,\ldots,x_n\in [-1,1]$ other than the Chebyshev extrema, such
that
$$
f[x_0,\ldots,x_n]\le 2^{n-1}||f||_{\infty},
$$
for all $f\in C[-1,1].$ By considering a suitably defined $f$, we see that
$$
\sum_{j=0}^n\prod_{k:k\neq j} |x_k-x_j|^{-1}\le 2^{n-1}.
$$

 We have $|T_n(x_j)|<1$ for some $j,$ as $x_0,\ldots,x_n$ are not the Chebyshev extrema. Thus
$$
\sum_{j=0}^n(-1)^{j+n} n!\prod_{k:k\neq j} |x_k-x_j|^{-1}T(x_j)< n!2^{n-1}.
$$
On the other hand,
$$
\sum_{j=0}^n(-1)^{j+n} n!\prod_{k:k\neq j} |x_k-x_j|^{-1}T(x_j)= n!2^{n-1}
$$
by Lemma~\ref{mvt}, and we have a contradiction.
\end{proof}

\section{A complex mean value theorem}

The following theorem is a complex version of the second mean value theorem for integrals. It would be unwise to suggest that such a classical result is new, but it seems possible that it is not well known.

\begin{cmvt}\label{cmvt}
Suppose that $f:[a,b]\to \mathbb{R}$ is monotone, $g:[a,b]\to \mathbb{C}$ is continuous, and $I=\int_a^bf(x)g(x)\,dx=|I|e^{i\theta}.$
Then there exists a point $c\in [a,b]$ such that
$$
I=\left(f(a)\R\left(e^{-i\theta}\int_a^cg(x)\,dx\right)+f(b)\R\left(e^{-i\theta}\int_c^bg(x)\,dx\right)\right)e^{i\theta},
$$
so that
$$
|I|\le \left|f(a)\int_a^cg(x)\,dx\right|+\left|f(b)\int_c^bg(x)\,dx\right|.
$$
Furthermore if $f$ is of constant sign and $|f|$ is decreasing, then there exists a point $c\in [a,b]$ such that
$$
|I|=f(a)\R\left(e^{-i\theta}\int_a^cg(x)\,dx\right)\le \left|f(a)\int_a^cg(x)\,dx\right|.
$$ 
\end{cmvt}
\begin{proof}
Let $g(x)=w(x)+iv(x)$ and 
$
I=\int_a^bf(x)g(x)\,dx=|I|e^{i\theta},
$
so that
\begin{align*}
|I|=Ie^{-i\theta}&=\cos\theta\int_a^bf(x)w(x)\,dx +\sin\theta\int_a^bf(x)v(x)\,dx\\
&=\int_a^bf(x)\left(\cos\theta w(x)+\sin\theta v(x)\right)\,dx.
\end{align*}
Now as $f(x)$ is monotonic and $(\cos\theta w+\sin\theta v)(x)$ is continuous, there exists a point $c\in[a,b]$ such that
\begin{align*}
|I|&=f(a)\int_a^c\left(\cos\theta w+\sin\theta v\right)(x)\,dx+f(b)\int_c^b\left(\cos\theta w+\sin \theta v\right)(x)\,dx\\
&=f(a)\R\left(e^{-i\theta}\int_a^cg(x)\,dx\right)+f(b)\R\left(e^{-i\theta}\int_c^bg(x)\,dx\right),
\end{align*}
by the second mean value theorem for integrals (see for example \cite{bar}).

Finally if $f$ is of constant sign and $|f|$ is decreasing we define $h$ to be equal to $f$ on $[a,b)$, set $h(b)=0$ and apply the argument to $h$.
\end{proof}

 The following corollary is perhaps the cornerstone of the theory
of oscillatory integrals, and is proven, using different methods, in J.G. van der Corput \cite{van}, A. Zygmund \cite{zy}, and E.M Stein \cite{st} with constants $2\sqrt{2}$, 4,
and 3 respectively. We use Theorem~\ref{cmvt} to find
the sharp constant.

\begin{vd1}\label{vd1} Suppose that $f:[a,b]\to \mathbb{R}$ is differentiable on $(a,b)$, $f'$ is increasing, and   $f'\ge \la >0$  
on $[a,b]$. Then  
$$
|I|=\left|\int_a^be^{if(x)}\,dx\right| \le \frac{1+\sin(\theta-f(a))}{\la}\le \frac{2}{\la},
$$
where $I=|I|e^{i\theta}$.
\end{vd1}

\begin{proof} By a change of variables we have
$$
I=\int_a^be^{if(x)}\,dx=\int_{f(a)}^{f(b)}\frac{e^{iy}}{f'(f^{-1}(y))}\,dy.
$$
Now as $1/(f'\circ f^{-1})$ is positive and decreasing, there exists a point $c\in [a,b]$ such that 
$$
I=\frac{1}{f'(a)}\R\left(e^{-i\theta}\int_{f(a)}^ce^{iy}\,dy\right)=\frac{1}{f'(a)}(\sin(c-\theta)-\sin(f(a)-\theta)),
$$
by Theorem~\ref{cmvt}, where $\theta$ is the direction of the integral. Thus
$$
|I|\le\frac{1+\sin(\theta-f(a))}{\la},
$$
as desired.
\end{proof}

We note that if the integral is non-zero, then its argument cannot be $f(a)+3\pi/2$. That the constant 2 is sharp is observed by considering $f(x)=x$ on $(0,\pi)$ with $\la=1$ . 
Similarly we obtain a Reimann--Lebesgue lemma.

\begin{reim} Suppose that $f\in L^1[0,1]$ is increasing. Then
$$
|\hat{f}(n)|=\left|\int_0^1f(x)e^{-2\pi i nx}\,dx\right|\le \frac{f(1)-f(0)}{2\pi n}(1-\sin \theta),
$$
where $\hat{f}(n)=|\hat{f}(n)|e^{i\theta }.$
\end{reim}

\begin{proof} Let $\hat{f}(n)=|\hat{f}(n)|e^{i\theta}.$ By a change of variables,
$$
\hat{f}(n)=\frac{1}{2\pi n}\int_0^{2\pi n}f\left(\frac{y}{2\pi n}\right)e^{- i y}\,dy,
$$
so by Theorem~\ref{cmvt}, there exists a point $c\in [0,2\pi n]$ such that
\begin{align*}
|\hat{f}(n)|&=\frac{f(0)}{2\pi n}\R\left(e^{-i\theta}\int_0^ce^{-iy}\,dy\right)+\frac{f(1)}{2\pi n}\R\left(e^{-i\theta}\int_c^{2\pi n}e^{-iy}\,dy\right)\\
&=\frac{f(0)}{2\pi n}(\sin \theta - \sin(c+\theta))+ \frac{f(1)}{2\pi n}(\sin(c+\theta)-\sin \theta).
\end{align*}
Thus
$$
|\hat{f}(n)|\le \frac{f(1)-f(0)}{2\pi n}(1-\sin \theta),
$$
as desired.
\end{proof}
Again we note that the Fourier transforms of an increasing function cannot have argument $3\pi/2.$ 

It is evident that there are versions of the above results for both increasing and decreasing functions.
%
%
%
%
%
\section{Van der Corput lemmas}\label{threet}

The following lemma, with constants taken to be $$C_n=2^{5/2}\pi^{1/n}(1-1/n),$$ is due to G.I. Arhipov, A.A.
Karacuba and V.N. Cubarikov \cite{ar}. It is a more precise formulation of what is generally known as van der Corput's Lemma. We improve their constants, so that the bound becomes sharp as $n$ tends to infinity. 

\begin{vd}\label{vd} Suppose that $f:(a,b)\to \mathbb{R}$ is $n$ times differentiable, where $n\ge2,$ and $|f^{(n)}(x)|\ge \la >0$
on $(a,b)$. Then  
$$
\left|\int_a^be^{ if(x)}\,dx\right| \le C_n\frac{n}{\la^{1/n}},
$$
where $C_n\le 2^{5/3}$ for all $n\ge 2$ and $C_n\to 4/e$ as $n\to \infty.$
\end{vd}
\begin{proof}
We integrate over $E_1=\{ x\in (a,b) : |f'(x)|\le \alpha \}$ and $E_2=\{ x\in (a,b) : |f'(x)|> \alpha \}$
separately. If we consider $g=f',$ we have $|g^{(n-1)}|\ge \lambda,$ so that
$$
|E_1|\le \left((n-1)!2^{2(n-1)-1}\right)^{1/(n-1)}\left(\frac{\alpha}{\lambda}\right)^{1/(n-1)},$$
by Lemma~\ref{twot}. So trivially 
$$
\left|\int_{E_1}e^{ if(x)} \,dx\right|\le \left((n-1)!2^{2n-3}\right)^{1/(n-1)}\left(\frac{\alpha}{\lambda}\right)^{1/(n-1)}.$$
Now $E_2$ is made up of at most $2(n-1)$ intervals on each of which $|f'|\ge \alpha$ and $f'$ is monotone, so by
Lemma~\ref{vd1} we have  
$$
\left|\int_{E_2}e^{if(x)} \,dx\right|\le 2(n-1)\frac{2}{\alpha}.
$$
Finally we optimise with respect to $\alpha$ and deduce that
\begin{equation}\label{lem2}
\left|\int_a^be^{if(x)}\,dx\right| \le \left(\frac{(n-1)!2^{2n-1}}{(n-1)^{n-2}}\right)^{1/n}\frac{n}{\lambda^{1/n}}.
\end{equation}
This bound tends to $(4/e)n/\lambda^{1/n}$ by Stirling's formula and we can check the first few terms to see it is always
less than or equal to the stated bound.
\end{proof}

To see that the bound becomes sharp as $n$ tends to infinity, we consider $f_n$ defined on $[-1,1]$ by
$$
f_n(x)=\frac{T_n(x)}{n}.
$$
 When $n\ge 2,$ we have
$$
 \left|\int_{-1}^1e^{ if_n(x)}\,dx\right|\ge\left|\int_{-1}^1\cos (f_n(x))\,dx\right|\ge \left|\int_{-1}^11-\frac{(f_n(x))^2}{2}\,dx\right|,
$$
by Taylor's Theorem. Now as $|f_n|\le 1/n,$ we see that
$$
\left|\int_{-1}^1e^{ if_n(x)}\,dx\right|\ge 2-\frac{1}{n^2},
$$
which tends to 2 as $n$ tends to infinity. On the other hand $$f_n^{(n)}=(n-1)!2^{n-1},$$
so that
$$
\left|\int_{-1}^1e^{ if_n(x)}\,dx\right|\le\left(\frac{(n-1)!2^{2n-1}}{(n-1)^{n-2}}\right)^{1/n}\frac{n}{((n-1)!2^{n-1})^{1/n}},
$$
by (\ref{lem2}) in Lemma~\ref{vd}.
We manipulate to obtain
$$
\left|\int_{-1}^1e^{ if_n(x)}\,dx\right|\le\left(\frac{2^{n}n^n}{(n-1)^{n-2}}\right)^{1/n}.
$$
This bound also tends to 2 as $n$ tends to infinity. Hence Lemma~\ref{vd} is asymptotically sharp.

We note that as
$$
\left|\int_{-2}^2e^{ix^2/2}\,dx\right|=3.33346...>\frac{4}{e}\times2
$$
and
$$
\left|\int_{-3}^{3}e^{i(x^3/6-x)}\,dx\right|=4.61932...>\frac{4}{e}\times3,
$$
the asymptote could not come strictly from below.

\begin{cor}\label{cor} Let $f(x)=a_0+a_1x+\dots+a_nx^n$ be a real polynomial of degree $n\ge 1.$ Then for all $a,b\in \mathbb{R}$,
$$
\al \int_a^be^{if(x)}\,dx\ar < \frac{C_n}{|a_n|^{1/n}},
$$
 where $C_n<11/2$ for all $n\ge 1$ and $C_n\to 4$ as $n\to \infty.$
\end{cor}
\begin{proof}
When $n=1,$ we obtain the result simply by integrating. For higher degrees we have $f^{(n)}=n!a_n,$ so that
$$
|I|=\al \int_a^be^{if(x)}\,dx\ar \le \left(\frac{(n-1)!2^{2n-1}}{(n-1)^{n-2}}\right)^{1/n}\frac{n}{(n!|a_n|)^{1/n}}.
$$
by (\ref{lem2}) in the proof of
Lemma~\ref{vd}. We can manipulate to obtain
$$
|I|\le \left(\frac{2^{2n-1}n^{n-1}}{(n-1)^{n-2}}\right)^{1/n}\frac{1}{|a_n|^{1/n}}.
$$
This bound tends to $4/|a_n|^{1/n}$ and we can check the first few terms to see it is
always less than the stated bound.
\end{proof}

Finally we improve the bound in Lemma~\ref{vd} when $n=2.$ This was the lemma that van der Corput originally proved with constant $2^{7/4}2\approx6.73.$ 

\begin{enequals2} Suppose that $f:(a,b)\to \mathbb{R}$ is twice differentiable and $|f^{''}(x)|\ge \la >0$
on $(a,b)$. Then  
$$
\left|\int_a^be^{ if(x)}\,dx\right| \le \frac{3^{3/4}2}{\sqrt{\la}}.
$$
\end{enequals2}
\begin{proof}
Let $I=\int_a^be^{ if(x)}\,dx=|I|e^{i\theta},$ so that
$$
|I|=e^{-i\theta}\int_a^be^{ if(x)}\,dx=\int_a^b\cos(f(x)-\theta)\,dx=\int_{a-\theta}^{b-\theta}\cos f(x)\,dx.
$$
We have reduced the problem to a real one, and from now on we shall consider
$$
I=\int_{a}^{b}\cos f(x)\,dx.
$$

Suppose that there exists a point $\zeta\in (a,b)$ such that $f'(\zeta)=0.$  (The proof when there is no such point is easier as we only need split the integral into two pieces
rather than three.) We suppose, without loss of generality, that
$$
\int_{a}^{\zeta}\cos f(x)\,dx\ge\int_{\zeta}^{b}\cos f(x)\,dx,
$$
and define $g$ by
$$
g(x)=\left\{ \begin{array}{cc}
f(x) & a\le x\le \zeta\\
f(2\zeta-x) & d\le x \le 2\zeta-a.
\end{array} \right.
$$
Now $|g''|>\la$ on $(a,2\zeta-a),$ and 
$$
\int_{a}^{b}\cos f(x)\,dx\le\int_{a}^{2\zeta-a}\cos g(x)\,dx,
$$
so from now on we will suppose that $f$ is even around $\zeta$. We will also suppose, without loss of generality, that $f''<-
\la.$

We will integrate over $(a,c)=\{x\in(a,\zeta):|f'(x)|>\alpha\},$ $(c,d)=\{x\in(a,b):|f'(x)|\le\alpha\},$ and $(d,b)=\{x\in(\zeta,b):|f'(x)|>\alpha\}$ separately. We have 
\begin{equation*}
\int_{a}^{c}\cos f(x)\,dx=\int_{f(a)}^{f(c)}\frac{\cos y}{f'(f^{-1}(y)}\,dy\le\frac{1+\sin f(c)}{\alpha},
\end{equation*}
by the second mean value theorem, as in Corollary~\ref{vd1}. Similarly,
$$
\int_{d}^{b}\cos f(x)\,dx\le \frac{1+\sin f(d)}{\alpha}=\frac{1+\sin f(c)}{\alpha},
$$
as $f$ is even around $\zeta$. By Lemma~\ref{twot}, we have 
$$
\int_c^d\cos f(x)\,dx\le \frac{2\alpha}{\la}\max_{x\in(c,d)} \cos f(x),
$$
and as $|f|\le \alpha$ in $(c,d)$ we see that
$$
|I|\le  \frac{2(1+\sin f(c))}{\alpha}+\frac{2\alpha}{\la}\max_{y\in(f(c),\alpha^2/\la)} \cos y.
$$
If we consider $\sqrt\la\le \alpha \le \sqrt{2\la},$ it is clear that this expression is maximised when $f(c)\in [0,\pi/2]$, and that $\max_{y\in(f(c),\alpha^2/\la)} \cos y=\cos f(c).$ Thus
\begin{equation*}\label{due}
|I|\le \max_{\theta\in [0,\pi/2]} \frac{2(1+\sin \theta)}{\alpha} +\frac{2\alpha}{\la}\cos \theta.
\end{equation*}
We optimise with respect to $\alpha$, so that
$$\alpha =\sqrt{\frac{1+\sin \theta}{\cos \theta}}\sqrt{\la}.
$$ Thus
$$
|I|\le \frac{2}{\sqrt{\la}}\thinspace\max_{\theta\in[0,\pi/2]}\sqrt{(\cos \theta, \sin \theta)\cdot(1,\sin \theta)}\le \frac{3^{3/4}2}{\sqrt{\la}},
$$
where the maximum is attained when $\theta=\pi/6.$
\end{proof}

In seems unlikely that the constant $3^{3/4}2\approx 4.56$ is sharp. Indeed some numerical experiments have suggested that the sharp constant when $n=2$ is 
$$2\sqrt{2}\max_{\theta\in[0,2\pi]}(\cos\theta,\sin\theta)\cdot\left(\int_0^{\sqrt{\pi/2+\theta}}\cos x^2\,dx,\int_0^{\sqrt{\pi/2+\theta}}\sin x^2\,dx\right),
$$
which is approximately equal to $3.3643....$
Hence any expression for the sharp constant would provide information about the Fresnel functions. 
 
 When $n=3$, the polynomials
$a_1 x+a_3x^3$ with the property $a_3/a_1^3=-0.3547...$ appear to be optimal for maximising 
\begin{equation}\label{pop}
\max_{a,b}\left|\int_{a}^be^{i(a_1x+a_3x^3)}\,dx\right|a_3^{1/3}.
\end{equation}
They share the property that they have local maxima that take the values\\ $\pm0.5935...$. These polynomials appear to be non-standard, so it may be difficult to find the absolutely
sharp constants in Lemma~\ref{vd} or Corollary~\ref{cor}. On the other hand, as the maximum in (\ref{pop}) appears to be $2.6396...<4$, it may be possible to show that the asymptote in
Corollary~\ref{cor} comes from below. 

\ \newline
\noindent Thanks to M. Cowling for all his invaluable help.


\begin{thebibliography}{9}


\bibitem{ar} G.I. Arhipov, A.A. Karacuba and V.N. Cubarikov, Trigonometric integrals, \textit{Math. USSR Izvestija} \textbf{15} (1980), 211--239. 
\bibitem{bar} R. G. Bartle, {\it The elements of real analysis}, Second edition, John Wiley \& Sons, New York, 1976. MR0393369 
\bibitem{ca} A. Carbery, M. Christ\ and\ J. Wright, Multidimensional van der Corput and sublevel set estimates, J. Amer. Math. Soc. {\bf 12} (1999), no.~4, 981--1015. MR1683156 
\bibitem{cw} A. Carbery\ and\ J. Wright, What is van der Corput's lemma in higher dimensions?, Publ. Mat. {\bf 2002}, Vol. Extra, 13--26. MR1964813
\bibitem{van} J.G. van der Corput, Zahlentheoretische absch\"atzungen, \textit{Math. Ann.} \textbf{84} (1921), 53--79.
\bibitem{ik} E. Isaacson\ and\ H. B. Keller, {\it Analysis of numerical methods}, Wiley, New York, 1966. MR0201039 
\bibitem{th} T. J. Rivlin, {\it Chebyshev polynomials}, Second edition, Wiley, New York, 1990. MR1060735 
\bibitem{st} E. M. Stein, {\it Harmonic analysis: real-variable methods, orthogonality, and oscillatory integrals}, Princeton Univ. Press, Princeton, NJ, 1993. MR1232192 
\bibitem{zy} A. Zygmund, {\it Trigonometric series. 2nd ed. Vols. I, II}, Cambridge Univ. Press, New York, 1959. MR0107776 












\end{thebibliography}
\end{document}